\documentclass[11pt]{article}
\usepackage{cite}
\usepackage{mathrsfs}
\usepackage{amssymb,amsfonts,amsmath}
\usepackage{enumerate}
\usepackage{indentfirst}
\usepackage{latexsym,bm}
\usepackage[noend]{algpseudocode}
\usepackage{algorithmicx,algorithm}
\usepackage{algorithm}
\usepackage{makeidx}
\usepackage{fancybox}
\usepackage{color}
\usepackage{multicol}
\usepackage[latin1]{inputenc}
\usepackage{graphicx}
\usepackage{pstricks,pst-node,pst-text,pst-3d}
\usepackage{tikz}
\newtheorem{thm}{Theorem}[section]
\newtheorem{cor}[thm]{Corollary}
\newtheorem{lem}[thm]{Lemma}
\newtheorem{op}[thm]{Problem}

\newtheorem{conj}{Conjecture}[section]

\newenvironment{pf}{{\noindent \it \bf Proof:}}{{\hfill$\Box$}\\}

\begin{document}

\title{\bf  Packing Strong Subgraph in Digraphs}
\author{
\small Yuefang Sun$^{1}$, Gregory Gutin$^{2}$, Xiaoyan Zhang$^{3}$\footnote{Corresponding author.
\vskip1mm Email address:
sunyuefang@nbu.edu.cn (Sun), gutin@cs.rhul.ac.uk (Gutin), zhangxiaoyan@njnu.edu.cn (Zhang)} \\
\small $^{1}$ School of Mathematics and Statistics,\\
\small  Ningbo University, Ningbo 315211, P. R. China \\
\small $^{2}$ Department of Computer Science, Royal Holloway,\\
\small University of London, Egham, UK\\
\small $^{3}$ School of Mathematical Science \& Institute of Mathematics,\\
\small Nanjing Normal University, Nanjing 210023, P. R. China\\
}


\date{}
\maketitle

\begin{abstract}

In this paper, we study two types of strong subgraph packing problems in digraphs, including internally disjoint strong subgraph packing problem and arc-disjoint strong subgraph packing problem. These problems can be viewed as generalizations of the famous Steiner tree packing problem and are closely related to the strong arc decomposition problem. We first prove the NP-completeness for the internally disjoint strong subgraph packing problem restricted to symmetric digraphs and Eulerian digraphs. Then we get inapproximability results for the arc-disjoint strong subgraph packing problem and the internally disjoint strong subgraph packing problem. Finally we study the arc-disjoint strong subgraph packing problem restricted to digraph compositions and obtain some algorithmic results by utilizing the structural properties.
\vspace{0.3cm}
\vspace{0.3cm}
\end{abstract}

\section{Introduction}
We refer the readers to \cite{Bang-Jensen-Gutin} for graph-theoretical notation and terminology not given
here. The Steiner Type Problems have attracted significant attention from researchers due to their importance in theoretical research and practical implications\cite{Henning-Nielsen-Oellermann,Grotschel-Martin-Weismantel,Sherwani,Li-Mao5,Bang-Jensen-Gutin-Yeo1,West-Wu,Bresar2011}. For a graph $G=(V,E)$ and a set $S\subseteq V$ of at least two
vertices, an {\em $S$-Steiner tree} is a tree $T$ of $G$ with $S\subseteq V(T)$. The basic problem of {\sc
Steiner Tree Packing} is to find a largest collection of edge-disjoint
$S$-Steiner trees in a given undirected graph $G$. Besides this classical version, some its variations were also studied, see e.g.
\cite{Cheriyan-Salavatipour, DeVos-McDonald-Pivotto, Kriesell, Lau,West-Wu}. The Steiner tree packing problem has applications in a number of areas such as VLSI circuit design \cite{Grotschel-Martin-Weismantel, Sherwani} and stream broadcasting \cite{Li-Mao5}. 

It is natural to consider 
extensions of the Steiner tree packing problem to directed graphs. One approach is to replace undirected tree with out-trees i.e. trees oriented from their roots, see e.g.
Cheriyan and Salavatipour \cite{Cheriyan-Salavatipour}, Sun and Yeo \cite{Sun-Yeo}.
In this paper, we will further study the {\sc
Strong Subgraph Packing} problem which can be considered as another extension of the Steiner tree packing problem and is closely related to the strong arc decomposition problem \cite{Bang-Jensen-Gutin-Yeo, Bang-Jensen-Huang, Bang-Jensen-Yeo, Sun-Gutin-Ai}. A digraph $D$ is {\em strong connected} (or, {\em strong}), if for any pair of vertices $x, y\in V(D)$, there is a path from $x$ to $y$ in $D$, and vice versa. Let $D=(V(D),A(D))$ be a digraph of order $n$, $S\subseteq V$ a
$k$-subset of $V(D)$ and $2\le k\leq n$. A strong subgraph $H$ of $D$ is called an {\em $S$-strong
subgraph} if $S\subseteq V(H)$. Two $S$-strong subgraphs are said to be {\em arc-disjoint} if they have no common arc. Furthermore, two
arc-disjoint $S$-strong subgraphs are said {\em internally disjoint}
if the set of common vertices of them is exactly $S$.

In this paper, we consider the following two types of strong subgraph packing problems in digraphs. The input of
{\sc Arc-disjoint strong subgraph packing} (ASSP)
consists of a digraph $D$ and a subset of vertices $S\subseteq
V(D)$, the goal is to find a largest collection of arc-disjoint
$S$-strong subgraphs. Similarly, the input of {\sc Internally-disjoint strong subgraph packing} (ISSP)
consists of a digraph $D$ and a subset of vertices $S\subseteq
V(D)$, and the goal is to find a largest collection of internally
disjoint $S$-strong subgraphs.


Let {\em internally (resp. arc-)disjoint strong subgraph packing number},
denoted by $\kappa_S(D)$~(resp. $\lambda_S(D)$), be the maximum number of internally
(resp. arc-)disjoint $S$-strong subgraphs in $D$. Then the {\em strong subgraph
$k$-connectivity} introduced in \cite{Sun-Gutin-Yeo-Zhang} is defined as
$$\kappa_k(D)=\min\{\kappa_S(D)\mid S\subseteq V(D), |S|=k\}.$$
Similarly, the {\em strong subgraph
$k$-arc-connectivity} introduced in \cite{Sun-Gutin2} is defined as
$$\lambda_k(D)=\min\{\lambda_S(D)\mid S\subseteq V(D), |S|=k\}.$$

A digraph $D$ is {\em symmetric} if it can be obtained from an undirected graph $G$ by replacing every edge $\{x,y\}$ of $G$ by the pair $xy,yx$ of edges. We denote it by $D=\overleftrightarrow{G}.$
A digraph $D$ is {\em Eulerian} if its undirected underlying graph is connected and the out-degree and in-degree of each vertex of $D$ coincide.
Clearly, symmetric digraphs is a subfamily of Eulerian digraphs.

It is worth mentioning that strong subgraph
connectivity is related to other concepts in graph theory. It is an extension of the well-established
tree connectivity of undirected graphs \cite{Li-Mao5}. Since $\kappa_2(\overleftrightarrow{G})=\kappa(G)$ \cite{Sun-Gutin-Yeo-Zhang} and $\lambda_2(\overleftrightarrow{G})=\lambda(G)$ \cite{Sun-Gutin2}, $\kappa_k(D)$ and
$\lambda_k(D)$ could be seen as generalizations of connectivity and
edge-connectivity of undirected graphs, respectively. 
The following concept of strong arc decomposition (or good decomposition) was studied in \cite{Bang-Jensen-Gutin-Yeo, Bang-Jensen-Huang, Bang-Jensen-Yeo, Sun-Gutin-Ai}. A digraph $D=(V,A)$ has a {\em strong arc decomposition} if $A$ has two disjoint sets $A_1$ and $A_2$ such that both $(V,A_1)$ and $(V,A_2)$ are strong. By definition, $\kappa_n(D)\geq 2$ (or, $\lambda_n(D)\geq 2$) if and only if $D$ has a strong arc decomposition.

In \cite{Sun-Gutin-Yeo-Zhang, Sun-Gutin2}, some hardness results for the decision versions of {\sc ISSP} and {\sc ASSP} were obtained. For general digraphs, we list such results in the following two tables.

\begin{center}
\begin{tabular}{|c|c|c|} \hline
\multicolumn{3}{|c|}{General digraphs} \\ \hline
$\lambda_{S}(D) \geq \ell$?                                     & $k \geq 2$     &  $k$ part \\
$ |S|=k$                                                              & constant       &  of input \\ \hline  \hline

$\ell \geq 2$ constant                                     & NP-complete \cite{Sun-Gutin-Yeo-Zhang}   &  NP-complete\cite{Sun-Gutin-Yeo-Zhang}  \\ \hline
$\ell$ part of input                                       & NP-complete\cite{Sun-Gutin-Yeo-Zhang}    &  NP-complete\cite{Sun-Gutin-Yeo-Zhang}  \\ \hline
\end{tabular}
\end{center}

\begin{center}
\begin{tabular}{|c|c|c|} \hline
\multicolumn{3}{|c|}{General digraphs} \\ \hline
$\kappa_{S}(D) \geq \ell$?                                     & $k \geq 2$     &  $k$ part \\
$ |S|=k$                                                              & constant       &  of input \\ \hline  \hline

$\ell \geq 2$ constant                                     & NP-complete \cite{Sun-Gutin2}   &  NP-complete\cite{Sun-Gutin2}  \\ \hline
$\ell$ part of input                                       & NP-complete\cite{Sun-Gutin2}    &  NP-complete\cite{Sun-Gutin2}  \\ \hline
\end{tabular}
\end{center}



A digraph $D$ is {\em semicomplete} if at least one of the arcs $xy,yx$ is in $D$ for every distinct $x,y\in
V(D)$. A digraph $D$ is {\em locally in-semicomplete} if all in-neighbours of each vertex in $D$ induce a semicomplete digraph. Hence, a semicomplete digraph is also locally in-semicomplete.
A digraph $D$ is {\em symmetric} if there
is an opposite arc $yx$ for every arc $xy$. 
For semicomplete digraphs and symmetric digraphs, the following hardness results were obtained in \cite{Sun-Gutin-Yeo-Zhang}.

\begin{thm}\label{thm03}\cite{Sun-Gutin-Yeo-Zhang}
Let $k,\ell \geq 2$ be fixed integers. Let $D$ be a semicomplete digraph and $S \subseteq V(D)$ with $|S|=k$. The problem of
deciding whether $\kappa_{S}(D) \geq \ell$ is polynomial-time solvable.
\end{thm}

\begin{thm}\label{thm02}\cite{Sun-Gutin-Yeo-Zhang}
Let $k\ge 3$ be a fixed integer. Let $D$ be a symmetric digraph and $S \subseteq V(D)$ with $|S|=k$. The problem of
deciding whether $\kappa_{S}(D) \geq \ell$ is NP-complete, where $\ell \geq 1$ is an integer.
\end{thm}

\begin{thm}\label{thm010}\cite{Sun-Gutin-Yeo-Zhang}
Let $k,\ell \geq 2$ be fixed integers. Let $D$ be a symmetric digraph and $S \subseteq V(D)$ with $|S|=k$. The problem of
deciding whether $\kappa_{S}(D) \geq \ell$ is polynomial-time solvable.
\end{thm}

Let $D$ be a digraph with $V(D)=\{u_i\mid 1\leq i\leq t\}$ and let
$H_1,\dots, H_t$ be digraphs with $V(H_i)=\{u_{i,j_i}\mid 1\le j_i\le n_i\}$. In the rest of this paper, we set $n_0=\min\{n_i\mid
1\leq i\leq t\}$. The {\em composition} of $D$ and $H_i$,
denoted by $Q=D[H_1,\dots , H_t]$, is a digraph with vertex set $\bigcup_{i=1}^t V(H_i)=\{u_{i,j_i}\mid 1\le i\le t, 1\le j_i\le n_i\}$ and arc set
\[
\left(\bigcup^t_{i=1}A(H_i) \right) \bigcup  \left( \bigcup_{u_iu_p\in A(D)} \{u_{i,j_i}u_{p,q_p} \mid 1\le j_i\le n_i, 1\le q_p\le n_p\} \right).
\]
A digraph $D$ is {\em quasi-transitive}, if for any triple $x,y,z$ of
distinct vertices of $D$, the following holds: if $xy$ and $yz$ are arcs of $D$ then
either $xz$ or $zx$ or both are arcs of $D$. Based on the notion of digraph composition, Bang-Jensen and Huang gave a recursive characterization of quasi-transitive digraphs below, where the decomposition is called the {\em canonical decomposition} of a quasi-transitive digraph.

\begin{thm}\label{intro01}\cite{Bang-Jensen-Huang2}
Let $D$ be a quasi-transitive digraph. Then the following assertions hold:
\begin{description}
\item[(a)] If $D$ is not strong, then there exists a transitive oriented
graph $T$ with vertices $\{u_i\mid i\in [t]\}$ and strong quasi-transitive digraphs $H_1, H_2, \dots, H_t$ such that $D = T[H_1, H_2, \dots, H_t]$, where $H_i$ is substituted for $u_i$, $i\in [t]$.
\item[(b)] If $D$ is strong, then there exists a strong semicomplete
digraph $S$ with vertices $\{v_j\mid j\in [s]\}$ and quasi-transitive digraphs $Q_1, Q_2, \dots, Q_s$ such that $Q_j$ is either a vertex or is non-strong and $D = S[Q_1, Q_2, \dots, Q_s]$, where $Q_j$ is substituted for $v_j$, $j\in [s]$.
\end{description}
\end{thm}

Composition of digraphs is a useful concept and tool in digraph theory, especially in 
the structural and algorithmic applications for quasi-transitive digraphs and their extensions, see
e.g. \cite{Bang-Jensen-Gutin, Bang-Jensen-Huang2, GSHC}.
In addition, digraph compositions generalize some families of digraphs, including (extended) semicomplete digraphs, quasi-transitive digraphs (by Theorem~\ref{intro01}) and lexicographic product digraphs (when $H_i$ is the same digraph $H$ for every $i\in [t]$, $Q$ is the lexicographic product of $T$ and $H$, see, e.g., \cite{Hammack}). In particular, semicomplete compositions generalize strong quasi-transitive digraphs. To see that strong compositions form a significant generalization of strong quasi-transitive digraphs, observe that the Hamiltonian cycle problem is polynomial-time solvable for quasi-transitive digraphs \cite{Gutin4}, but NP-complete for strong compositions (see, e.g., \cite{Bang-Jensen-Gutin-Yeo}). There are several papers appeared on the topic of digraph compositions \cite{Bang-Jensen-Gutin-Yeo, Gutin-Sun, Sun-Gutin-Ai}. 

The rest of the paper is organized as follows. In Section 2, we study the decision versions of {\sc ISSP} for symmetric digraphs (Theorem~\ref{thmSym}) and Eulerian digraphs (Theorem~\ref{thmEul}).
The new results together  with Theorems~\ref{thm02} and~\ref{thm010} allow us to
complete the following two tables, which clearly demonstrate that {\sc ISSP} increases in hardness when moving from symmetric to Eulerian digraphs.

\begin{center}
\begin{tabular}{|c||c|c|c|} \hline
\multicolumn{4}{|c|}{Table 1: Symmetric digraphs} \\ \hline
$\kappa_{S}(D) \geq \ell$? & $k=2$              & $k \geq 3$     &  $k$ part \\
$ |S|=k$                   &                    & constant       &  of input \\ \hline  \hline
$\ell \geq 2$ constant     & Polynomial \cite{Sun-Gutin-Yeo-Zhang}        & Polynomial  \cite{Sun-Gutin-Yeo-Zhang}   &  NP-complete  \\ \hline
$\ell$ part of input       & Polynomial \cite{Sun-Gutin-Yeo-Zhang}         & NP-complete \cite{Sun-Gutin-Yeo-Zhang}   &  NP-complete \cite{Sun-Gutin-Yeo-Zhang} \\ \hline
\end{tabular}
\end{center}

\begin{center}
\begin{tabular}{|c||c|c|c|} \hline
\multicolumn{4}{|c|}{Table 2: Eulerian digraphs} \\ \hline
$\kappa_{S}(D) \geq \ell$? & $k=2$              & $k \geq 3$     &  $k$ part \\
$ |S|=k$                   &                    & constant       &  of input \\ \hline  \hline
$\ell \geq 2$ constant     & NP-complete        & NP-complete   &  NP-complete  \\ \hline
$\ell$ part of input       & NP-complete          & NP-complete   &  NP-complete  \\ \hline
\end{tabular}
\end{center}
In Section~3, we obtain inapproximability results on {\sc ISSP} and {\sc ASSP} in Theorem \ref{thma}.
Some structural properties and algorithmic results on digraph compositions will be given in Theorems \ref{thmc} in Section~4. We also discuss the min-max relation on strong subgraph packing problem and pose an open problem.

\section{Complexity for $\kappa_{S}(D)$ on symmetric digraphs and Eulerian digraphs}

\begin{thm}\label{thmSym}
Let $\ell \geq 2$ be a fixed integer. Let $D$ be a symmetric digraph and $S \subseteq V(D)$ ($k=|S|$ is part of the input). The problem of
deciding whether $\kappa_{S}(D) \geq \ell$ is NP-complete.
\end{thm}
\begin{pf} It is easy to see that this problem is in NP.
We will reduce from the NP-complete problem of $2$-coloring of hypergraphs (see \cite{Lovasz}).
That is, we are given a hypergraph $H$ with vertex set $V(H)$ and edge set $E(H)$, and want to
determine if we can $2$-colour the vertices $V(H)$ such that every hyperedge in $E(H)$ contains
vertices of both colours. 

Define a symmetric digraph $D$ as follows.
Let $U = \{u_1,u_2,\ldots, u_{\ell-2} \}$ and
let $V(D)= V(H) \cup E(H) \cup U \cup \{r\}$ and let the arc set of $D$ be defined as follows.

\[
\begin{array}{rcl}
A(D) & = & \{xe ,ex \; | \; x\in V(H), \mbox{ } e \in e(H) \mbox{ and } x \in V(e)\} \\
     &   &  \cup \; \; \{r u_i, u_i r, u_i e, e u_i \; | \; u_i \in U \mbox{ and } e \in E(H) \} \\
     &   &  \cup \; \;  \{r x, x r \; | \; x \in V(H)\} \\
\end{array}
\]

Let $S = E(H) \cup \{r\}$. This completes the construction of $D$ and $S$.
We will show  that $\kappa_{S}(D)\geq \ell$ if and only if $H$ is 2-colourable, which
will complete the proof.

First assume that $H$ is $2$-colourable and let $R$ be the red vertices in $H$ and $B$ be the blue vertices in $H$ in
a proper $2$-colouring of $H$. For $i=1,2,\ldots, \ell-2$, let $D_i$ contain all arcs between $u_i$ and $S$.
Let $D_{\ell-1}$ contain all arcs between $r$ and $R$, and for each edge $e \in E(H)$ we add all arcs between $R$ and $e$ in $D$ to $D_{\ell-1}$
(this is possible as every edge in $H$ contains a red vertex). Analogously, let $D_{\ell}$ contain all arcs between $r$ and $B$, and for each edge $e \in E(H)$ we add all arcs between $B$ and $e$ in $D$ to $D_{\ell}$ (again, this is possible as every edge in $H$ contains a blue vertex). Observe that $D_1, D_2, \ldots, D_{\ell}$ are internally disjoint $S$-strong subgraphs in $D$, so $\kappa_{S}(D)\geq \ell$.

Conversely, assume that $\kappa_{S}(D)\geq \ell$ and let $D_1', D_2', \ldots, D_{\ell}'$ be a set of $\ell$ internally disjoint $S$-strong subgraphs in $D$.
At least two of these subgraphs contain no vertex from $U$ (as $|U| = \ell-2$). Without loss of generality assume that
$D_1'$ and $D_2'$ do not contain any vertex from $U$. Let $B'=V(D_1')\cap V(H)$ and $R'=V(D_2')\cap V(H)$. Observe that $B'\cap R'=\emptyset$ as $D_1'$ and $D_2'$ are internally disjoint. Every $e \in E(H)$ has at least one neighbour in $B'$ and $R'$, respectively. Therefore $H$ is $2$-colourable (any vertex in $H$ that is not in either $R'$ or $B'$ can be assigned arbitrarily to either $B'$ or $R'$). This completes the proof.
\end{pf}

Recall that it was proved in \cite{Sun-Gutin-Yeo-Zhang} that $\kappa_2(\overleftrightarrow{G})=\kappa(G)$, which means that $\kappa_2(\overleftrightarrow{G})$ can be computed in polynomial time. In fact, the argument also means that $\kappa_{\{x, y\}}(\overleftrightarrow{G})=\kappa_{\{x,y\}}(G)$, that is, the maximum number of disjoint $x-y$ paths in $G$, therefore can be computed in polynomial time. Then combining with Theorems~\ref{thm03}, \ref{thm02} and~\ref{thmSym}, we can complete all the entries of Table~1.

Sun and Yeo proved the NP-completeness of the {\sc Directed 2-Linkage} problem for Eulerian digraphs. {\sc Directed 2-Linkage} is an NP-hard problem which can be stated as follows: Given a digraph $D$ and four vertices $s_1,t_1,s_2,t_2$, decide whether there are vertex-disjoint paths from $s_1$ and $t_1$ and from $s_2$ to $t_2$.

\begin{thm}\label{link-eulerian}\cite{Sun-Yeo}
The 2-linkage problem restricted to Eulerian digraphs is NP-complete.
\end{thm}

Using Theorem~\ref{link-eulerian}, we will now prove the following result for Eulerian digraphs which gives Table~2.

\begin{thm}\label{thmEul}
Let $k,\ell \geq 2$ be fixed. Let $D$ be an Eulerian digraph and $S \subseteq V(D)$ with $|S|=k$. The problem of
deciding whether $\kappa_{S}(D) \geq \ell$ is NP-complete.
\end{thm}
\begin{pf} 
Let $(D,s_1,t_1,s_2,t_2)$ be an instance of  {\sc Directed 2-Linkage} restricted to Eulerian digraphs. Let us first construct a new digraph $D'$ by adding to $D$ vertices $x, y, r_1, r_2$ and arcs $$t_1x,xs_1, t_2y,ys_2,  xs_2,s_2x,yt_1,t_1y, s_1r_1, r_1t_2, s_2r_2, r_2t_1.$$

Secondly, we add to $D'$ $\ell -2$ copies of the 2-cycle $xyx$ and subdivide the arcs of every copy to avoid parallel arcs, that is, we insert each arc $xy$ (resp. $yx$) a new vertex $z_i$ (resp. $z'_i$) where $i\in [\ell-2]$. Let us denote the new digraph by $D''$. Note that $D''=D'$ for $\ell =2$.

Finally, we add to $D''$ $k-2$ new vertices $x_1,\dots ,x_{k-2}$ and arcs of $\ell$ 2-cycles $xx_ix$ for each $i\in [k-2]$. Subdivide the new arcs to avoid parallel arcs, that is, we insert each arc $xx_i$ (resp. $x_ix$) a new vertex $x_{i, j}$ (resp. $x'_{i, j}$) where $j\in [\ell]$. Let us denote the new digraph by $D'''$. Observe that $D'''$ is Eulerian as $D$ is Eulerian.

Let $S=\{x, y, x_i\mid i\in [k-2]\}$, $U=\{x_i, x_{i, j}, x'_{i, j}\mid i\in [k-2], j\in [\ell]\}$, $Z=\{z_j\mid j\in [\ell-2]\}$ and $Z'=\{z'_j\mid j\in [\ell-2]\}$.
It remains to show that $(D,s_1,t_1,s_2,t_2)$ is a positive instance of {\sc Directed 2-Linkage} restricted to Eulerian digraphs if and only if $\kappa_S(D''')\ge \ell$.

Suppose $(D,s_1,t_1,s_2,t_2)$ is a positive instance of {\sc Directed 2-Linkage} restricted to Eulerian digraphs, that is, there is a pair of vertex-disjoint $s_1-t_1$ path $P_1$ and $s_2-t_2$ path $P_2$. Let $H_1$ be the subdigraph of $D'''$ consisting of the arcs $xs_1, t_1x, t_1y, yt_1$, the path $P_1$, and the cycle $x, x_{i,\ell-1}, x_i, x'_{i,\ell-1}, x$ where $i\in [k-2]$. Let $H_2$ be the subdigraph of $D'''$ consisting of the arcs $ys_2, t_2y, s_2x, xs_2$ and the path $P_2$, and the cycle $x, x_{i,\ell}, x_i, x'_{i,\ell}, x$ where $1\leq i\leq k-2$. For $3\leq j\leq \ell$, let $H_j$ be the subdigraph of $D'''$ consisting of the cycles $x,z_{j-2},y,z'_{j-2},x$ and $x, x_{i,j-2}, x_i, x'_{i,j-2}, x$ where $i\in [k-2]$. Observe that $\{H_i\mid i\in [\ell]\}$ is a family of internally disjoint $S$-subgraphs, therefore $\kappa_S(D''')\ge \ell$.

If $\kappa_S(D''')\ge \ell$, then there is a set of internally disjoint $S$-subgraphs, say $\{H_i\mid i\in [\ell]\}$. Observe that $\{H'_i=H_i-U\mid i\in [\ell]\}$ is a set of $\ell$ internally disjoint $\{x,y\}$-subgraphs in $D''$.
Since $deg_{D''}^+(x)=deg_{D''}^-(x)=deg_{D''}^+(y)=deg_{D''}^-(y)=\ell$, each $H'_i$ contains precisely one out-neighbour and one in-neighbour of $x$ (resp. $y$). Therefore, there are two subdigraphs, say $H'_1, H'_2$, such that 
$V(H'_i)\cap Z=\emptyset$ for each $i\in [2]$. Let us consider two cases.

\vspace{2mm}

\noindent{\bf Case 1:} $V(H'_i)\cap Z'=\emptyset$ for each $i\in [2]$. Then we have  $V(H'_i)\subseteq V(D')$ for each $i\in [2]$.
Since the in-degree of $x$ in $D'$ is $2$, we may without loss of generality assume that
$t_1 \in V(H'_1)$ and $s_2 \in V(H'_2)$. As $y$ has in-degree $2$ in $D'$ and $t_1 \in V(H'_1)$ we must
have $t_2 \in V(H'_2)$. As the out-degree of $x$ is $2$, we analogously have $s_1 \in V(H'_1)$
(as $s_2 \in V(H'_2)$). Note that now we have  both $s_i$ and $t_i$ belong to $V(H'_i)$. Therefore, there must be a path $P_i$ from $s_i$ to $t_i$ in $H'_i$ and by definition of $D'$, $P_i$ will not have vertices outside of $D$. As $H_1$ and $H_2$ are internally disjoint, the paths are disjoint, so $(D,s_1,t_1,s_2,t_2)$ is a positive instance of {\sc Directed 2-Linkage} restricted to Eulerian digraphs.

\vspace{2mm}

\noindent{\bf Case 2:} $V(H'_i)\cap Z'\neq\emptyset$ for some $i\in [2]$. We will reduce Case 2 to Case 1. We just consider the case that $V(H'_1)\cap Z'\neq\emptyset$ and $V(H'_2)\cap Z'=\emptyset$ since the argument for the remaining case is similar. Let $V(H'_1)\cap Z'=\{z'_1\}$. Observe that there must exist some $H'_i$, say $H'_3$, such that $V(H'_3)\cap Z'=\emptyset$. Let $P'$ be a $y-x$ path in $H'_3$. Then we update $H'_1$ and $H'_3$ by exchanging the two paths $yz'_1x$ and $P'$ (note that in this procedure we may need to delete some vertices or arcs to guarantee the strong connectedness of updated $H'_1$ and $H'_3$ if necessary, and this will not affect the correctness), and we now also have $V(H'_i)\cap Z'=\emptyset$ for $i\in [2]$ and still make sure that the updated $\{H'_i\mid i\in [\ell]\}$ is a family of $\ell$ internally disjoint $\{x,y\}$-subgraphs in $D''$.
\end{pf}

Note that by Tables~1 and~2, for any fixed integers $k\geq 2$ and $\ell\geq 2$, the problem of deciding whether $\kappa_S(D)\geq \ell$ is NP-complete for an Eulerian digraph. However, when restricted to the class of symmetric digraphs, the above problem becomes polynomial-time solvable.



\section{Inapproximability results on {\sc ISSP} and {\sc ASSP}}
In the {\sc Set Cover Packing} problem, the input consists of a bipartite graph $G=(C\cup B, E)$, and the goal is to find a largest collection of pairwise disjoint set covers of $B$, where a {\em set cover} of $B$ is
a subset $S\subseteq C$ such that each vertex of $B$ has a neighbor
in $S$. Feige et al.
\cite{Feige-Halldorsson-Kortsarz-Srinivasan} proved the following inapproximability
result on the {\sc Set Cover Packing} problem.

\begin{thm}\label{thmFeige}\cite{Feige-Halldorsson-Kortsarz-Srinivasan}
Unless P=NP, there is no $o(\log{n})$-approximation algorithm for
{\sc Set Cover Packing}, where $n$ is the order of $G$.
\end{thm}

We now get our inapproximability results for {\sc ISSP} and {\sc ASSP} by reductions
from the {\sc Set Cover Packing} problem.

\begin{thm}\label{thma}The following assertions hold:\\
$(i)$~Unless P=NP, there is no $o(\log{n})$-approximation algorithm
for {\sc ISSP}, even restricted to the case that $D$ is a symmetric
digraph and $S$ is independent in $D$, where $n$ is the order of
$D$.\\
$(ii)$~Unless P=NP, there is no $o(\log{n})$-approximation algorithm
for {\sc ASSP}, even restricted to the case that $S$ is independent
in $D$, where $n$ is the order of $D$.
\end{thm}
\begin{pf}
\noindent{\bf Part (i)} Let $G(C\cup B, E)$ be an instance of {\sc
Set Cover Packing}. We construct an instance $(D, S)$ of {\sc ISSP}
by setting
$$V(D)=\{x\}\cup C\cup B,$$
$$A(D)=\{xu, ux\mid u\in C\}\cup \{uv, vu\mid
u~and~v~are~adjacent~in~G\}$$ and $$S=\{x\}\cup B.$$

If $\{C_i\subseteq C\mid 1\leq i\leq \ell\}$ is a set cover packing,
then the subdigraph in $D$ induced by the vertex set $\{x\}\cup C_i\cup B~(1\leq
i\leq \ell)$ forms a set of $\ell$ internally disjoint $S$-strong
subgraphs in $D$.

Conversely, let $\{D_i\mid 1\leq i\leq \ell\}$ be a set of $\ell$ internally
disjoint $S$-strong subgraphs in $D$. Since $B$ is an independent set
in $D$, for each $D_i$, there is a set $C_i\subseteq C$ of vertices satisfying the following:
every vertex in $B$ has a neighbor in $C_i$ such that it
can reach the vertex $x$. Observe that these sets $C_i$ are pairwise
disjoint and form a set cover packing of cardinality $\ell$. Note
that $D$ is symmetric and $S$ is an independent set of $D$. This
completes the proof of $(i)$ by Theorem \ref{thmFeige}.

\noindent{\bf Part (ii)} We construct an instance $(D', S')$ of {\sc
ASSP} from $(D, S)$ with $V(D')=\{x\}\cup B \cup \{u^+, u^-\mid u\in
C\}$ and
$S'=S=\{x\}\cup B$ such that:\\
$(1)$~$u^-u^+\in A(D')$ for each $u\in C$; \\
$(2)$~$vu^-\in A(D')$ if $vu\in A(D), v\in S', u\in
C$;\\
$(3)$~$u^+v\in A(D')$ if $uv\in A(D), v\in S', u\in C$.

If $\{C'_i\subseteq C\mid 1\leq i\leq \ell\}$ is a set cover
packing, then the subdigraph in $D'$ induced by the vertex set $\{x\}\cup
\{u^-,u^+\mid u\in C'_i\}\cup B~(1\leq i\leq \ell)$ forms a set of
$\ell$ arc-disjoint $S'$-strong subgraphs in $D'$.

Now let $\{D'_i\mid 1\leq i\leq \ell\}$ be a set of $\ell$
arc-disjoint $S'$-strong subgraphs in $D'$. In each $D'_i$, since
$B$ is an independent set in $D'$, each vertex $v\in B$ has to pass
through an arc of type $u^-u^+$ to reach $x$ for some $u\in C$.
Hence, in $G$ there is a set $C'_i\subseteq C$ of vertices such
that every vertex in $B$ has a neighbor in $C'_i$. Furthermore,
since the strong subgraphs $D'_i$ are pairwise arc-disjoint, the
sets $C'_i~(1\leq i\leq \ell)$ are pairwise disjoint and form a set
cover packing of cardinality $\ell$. Note that $S'$ is an
independent set of $D'$. This completes the proof of $(ii)$ by
Theorem \ref{thmFeige}.
\end{pf}

\section{Structural properties and algorithmic results on digraph compositions}

A digraph is {\em Hamiltonian decomposable}
if it has a family of Hamiltonian cycles such that every arc of
the digraph belongs to exactly one of the cycles. Ng \cite{Ng} proved the
following result on the Hamiltonian decomposition of complete
regular multipartite digraphs.

\begin{thm}\label{thm07}\cite{Ng}
The digraph $\overleftrightarrow{K}_{r, r, \ldots, r}$ ($s$ times)
is Hamiltonian decomposable if and only if $(r, s)\neq (4,
1)$ and $(r, s)\neq (6, 1)$.
\end{thm}

By Theorem \ref{thm07}, we can determine the precise value for the
strong subgraph $k$-arc-connectivity of a complete bipartite
digraph.

\begin{lem}\label{thm2}
For two positive integers $a$ and $b$ with $a\leq b$, we have
$$\lambda_k(\overleftrightarrow{K}_{a, b})=a$$ for $2\leq k\leq
a+b$.
\end{lem}
\begin{pf}
Let $V(\overleftrightarrow{K}_{a, b})= V_1\cup V_2$ with
$V_1=\{u_i\mid 1\leq i\leq a\}$ and $V_2=\{v_j\mid 1\leq j\leq b\}$.
By Theorem \ref{thm07}, the subgraph of
$\overleftrightarrow{K}_{a, b}$ induced by $\{u_i, v_j\mid 1\leq i,
j \leq a\}$ can be decomposed into $a$ Hamiltonian cycles:
$H_i~(1\leq i\leq a)$. For each $1\leq i\leq a$, let $D_i$ be the
strong spanning subgraph of $\overleftrightarrow{K}_{a, b}$
obtained from $H_i$ by adding the arc set $\{u_iv_j, v_ju_i\mid
a+1\leq j\leq b\}$. Observe that these subgraphs are
pairwise arc-disjoint, and so
$\lambda_{a+b}(\overleftrightarrow{K}_{a, b})\geq a$. It is known
\cite{Sun-Gutin2} that $\lambda_{k+1}(D)\leq \lambda_{k}(D)$ $(1\leq
k\le n-1)$ and $\lambda_k(D) \leq \min\{\delta^+(D), \delta^-(D)\}$
for a digraph $D$ with order $n$, we have that
$a=\min\{\delta^+(\overleftrightarrow{K}_{a, b}),
\delta^-(\overleftrightarrow{K}_{a, b})\}\geq
\lambda_{2}(\overleftrightarrow{K}_{a, b})\geq \ldots \geq
\lambda_{a+b}(\overleftrightarrow{K}_{a, b})\geq a$. This completes
the proof.
\end{pf}

The {\em lexicographic product} \cite{Hammack} $G\circ H$ of two digraphs $G$ and
$H$ is the digraph with vertex set $$V(G\circ H)=V(G)\times V(H)=\{(u,
u')\mid u\in V(G), u'\in V(H)\}$$ and arc set $$A(G\circ
H)=\{(u,u')(v,v')\mid uv\in A(G),~or~u=v~and~u'v'\in
A(H)\}.$$ The following
result was also proved by Ng, where $\overline{K_r}$ stands for the
digraph of order $r$ with no arcs and $\overrightarrow{C}_t$ is the directed cycle of order $t.$

\begin{lem}\label{thm08}\cite{Ng2}
For any two integers $t\geq 2$ and $r\geq 3$, the product digraph $\overrightarrow{C}_t\circ
\overline{K_r}$ is Hamiltonian decomposable.
\end{lem}

It follows from the constructive proof of  Lemma \ref{thm08} that the Hamiltonian cycles in the lemma can be found in $O(n^2)$ time. Recall that a strong semicomplete digraph is also locally in-semicomplete. By Camion's Theorem \cite{Camion}, there is a Hamiltonian cycle in a strong semicomplete digraph. In fact, Bang-Jensen and Hell obtained a stronger result.
\begin{thm}\label{thm12}\cite{Bang-Jensen-Hell}
There is an $O(m + n\log {n})$ algorithm for finding a
Hamiltonian cycle in a strong locally in-semicomplete digraph.
\end{thm}

Recall that $n_0=\min\{n_i\mid 1\leq i\leq t\}.$
By Lemma \ref{thm08} and Theorem \ref{thm12}, the following result holds.
\begin{lem}\label{thm4}
Let $Q=D[H_1,\dots ,H_t]$ with $|D|=t\geq 2$ and $|V(H_i)|\geq 3$ for each $1\leq i\leq t$. If $D$ is a strong semicomplete
digraph, then $Q$ has at least $n_0$ arc-disjoint strong
spanning subgraphs. Moreover, these strong subgraphs can be found in time
$O(n^2)$, where $n$ is the order of $Q$.
\end{lem}
\begin{pf}
By Theorem \ref{thm12}, we can find a Hamiltonian cycle of $D$ in time $O(n^2)$. Clearly,
$Q$ contains $\overrightarrow{C}_t\circ \overline{K_{n_0}}$ as
a spanning subgraph, where $t\geq 2$. By
Lemma \ref{thm08}, $\overrightarrow{C}_t\circ
\overline{K_{|V(H_1)|}}$ is Hamiltonian decomposible, and these Hamiltonian cycles can be found in time
$O(n^2)$. Furthermore, these cycles are desired strong spanning subgraphs in $Q$.
\end{pf}

Let $\mathcal{Q}_0=\{\overrightarrow{C}_3[\overline{K_2},\overline{K_2},\overline{K_2}], \overrightarrow{C}_3[\overrightarrow{P_2},\overline{K_2},\overline{K_2}], \overrightarrow{C}_3[\overline{K_2},\overline{K_2},\overline{K_3}]\}$. Sun, Gutin and Ai obtained the following characterization on arc-disjoint
strong spanning subgraphs in digraph compositions.
\begin{thm}\label{thm09}\cite{Sun-Gutin-Ai}
Let $D$ be a strong semicomplete digraph on $t\ge 2$ vertices and let $n_0\ge 2.$
Then $Q=D[H_1,\dots ,H_t]$ has a pair of arc-disjoint
strong spanning subgraphs if and only if $Q\not\in \mathcal{Q}_0$.
\end{thm}

We now give two sufficient conditions for a digraph composition to have at least $n_0$ arc-disjoint
$S$-strong subgraphs for any $S\subseteq V(Q)$ with $2\leq |S|\leq
|V(Q)|$.
\begin{thm}\label{thmc}
Let $Q=D[H_1,\dots, H_t]$ with $t\geq 2$. Then $Q$ has at least $n_0$ arc-disjoint
$S$-strong subgraphs for any $S\subseteq V(Q)$ with $2\leq |S|\leq
|V(Q)|$ if one of the following conditions holds:\\
$(i)$~$D$ is a strong symmetric digraph;\\
$(ii)$~$D$ is a strong semicomplete digraph and $Q\not\in \mathcal{Q}_0$.\\ 
Moreover, these strong subgraphs can be found in time
$O(n^4)$, where $n$ is the order of $Q$.
\end{thm}
\begin{pf}
\noindent{\bf Part (i)} For any $S\subseteq V(Q)$ with $2\leq |S|\leq
|V(Q)|$, we will obtain $n_0$ arc-disjoint $S$-strong subgraphs using the following three steps:

\noindent{\bf Step 1.} We obtain a spanning
subgraph $Q'=D[H'_1,\dots, H'_t]$ of $Q$ such that $V(H'_i)= V(H_i)$ and each $H'_i$ has no arcs, where $1\leq i\leq t$.

\noindent{\bf Step 2.}
For each pair of $1\leq p, q\leq t$ such that $u_pu_q, u_qu_p\in A(D)$, let
$Q_{p,q}$ be the subgraph of $Q$ induced by the vertex set
$\{u_{p, j_p}, u_{q, j_q} \mid 1\leq j_p\leq n_p, 1\leq j_q\leq n_q\}$. We obtain $n_0$ arc-disjoint strong spanning
subgraphs: $\{D_{p, q, s}\mid 1\leq s\leq n_0\}$ in $Q_{p,q}$ by the construction of Lemma~\ref{thm2}, since each $Q_{p,q}$ is a complete bipartite digraph.

\noindent{\bf Step 3.}
For each $1\leq
s\leq n_0$, let $D_s$ be the union of all $D_{p, q, s}$ with
$u_pu_q, u_qu_p\in A(D)$.

Observe that the subgraphs in Step 3
are strong and pairwise arc-disjoint, so we obtain a set of $n_0$
arc-disjoint strong spanning subgraphs of $Q'$, furthermore,
these subgraphs are desired arc-disjoint $S$-strong subgraphs of $Q$.

Step 1 can be performed in $O(n^2)$ time. In Step 2, there are at most $n\choose 2$ pairs of $p, q$, and note that $\{D_{p, q, s}\mid 1\leq s\leq n_0\}$ can be found in $O(n^2)$ time in $Q_{p,q}$ by the construction of Lemma
\ref{thm2}, so Step 2 can be executed in time $O(n^4)$. Step 3 can be performed in time $O(n^2)$.
Hence the desired
subgraphs can be found in polynomial time $O(n^4)$. This completes the
proof of part $(i)$.

\noindent{\bf Part (ii)} For the case that $n_0=1$, $Q$ itself is the desired strong subgraph.
The result holds for the case that
$n_0=2$ by Theorem \ref{thm09} and the fact that a strong spanning
subgraph is an $S$-strong subgraph for any $S\subseteq V(Q)$ with
$2\leq |S|\leq |V(Q)|$. It follows from the construction proof of Theorem \ref{thm09} that
these strong spanning subgraphs can be found in $O(n^3)$ time.

For the case that $n_0\geq 3$, we will get $n_0$ arc-disjoint $S$-strong subgraphs for any $S\subseteq V(Q)$ with $2\leq |S|\leq
|V(Q)|$ by the following two steps:

\noindent{\bf Step 1.} Find $n_0$
arc-disjoint strong spanning subgraphs: $D'_1, \ldots, D'_{n_0}$ in $Q'$ by Lemma \ref{thm4}, where $Q'=D[H'_1,\dots ,H'_t]$ is an
induced subgraph of $Q$ such that $V(H'_i)=\{u_{i,j_i}\mid 1\le
i\le t, 1\le j_i\le n_0\}$.

\noindent{\bf Step 2.}
For each $1\leq j\leq n_0$, we construct a spanning subgraph $D_j$
of $Q$ from $D'_j$ by adding arcs between $V(Q)\setminus V(Q')$ and
$\{u_{i,j}\mid 1\leq i\leq t\}$.

Observe that these subgraphs in Step 2 are strong and pairwise arc-disjoint, so we obtain a set of $n_0$ arc-disjoint strong spanning subgraphs of $Q$,
furthermore, these subgraphs are desired arc-disjoint $S$-strong subgraphs.
Step 1 can be performed in $O(n^2)$ time by Lemma \ref{thm08} and Step 2 can be performed in $O(n^3)$ time. This completes the proof of part $(ii)$.
\end{pf}

\begin{cor}\label{thm11}
Let $Q=D[H_1,\dots, H_t]$ with $t\geq 2$. Then $$\lambda_k(Q)\geq n_0$$ for any $2\leq k
\leq |V(Q)|$ if one of the following conditions holds:\\
$(i)$~$D$ is a strong symmetric digraph;\\
$(ii)$~$D$ is a strong semicomplete digraph and $Q\not\in \mathcal{Q}_0$.\\
Moreover, the bound is sharp in each case.
\end{cor}
\begin{pf} The bound clearly holds by Theorem \ref{thmc}.
For the sharpness of the bound for the first case, let $Q=D[\overline{K_r},\dots, \overline{K_r}]$ with
$|D|\geq 2$ and $D$ be a strong symmetric digraph with
$\min\{\delta^+(D), \delta^-(D)\}=1$. We clearly have
$\lambda_k(Q)\geq n_0=r$. Furthermore, by the fact that
$\lambda_k(Q)\leq \min\{\delta^+(Q), \delta^-(Q)\}=r$, we have
$\lambda_k(Q)=r$ for $2\leq k\leq |V(Q)|$.

For the sharpness of this bound for the second case, consider the
digraph $Q=\overrightarrow{C}_3[\overline{K_r},\dots,
\overline{K_r}]$ with $r\geq 3$. We clearly have $\lambda_k(Q)\geq
r$. Furthermore, by the fact that $\lambda_k(Q)\leq
\min\{\delta^+(Q), \delta^-(Q)\}=r$, we have $\lambda_k(Q)=r$ for
$2\leq k\leq |V(Q)|$.
\end{pf}

By definition, we have $G\circ H\cong G[H, \ldots, H]$. Then
by Lemma \ref{thm4}, Theorems \ref{thm09} and \ref{thmc}, the
following result directly holds:

\begin{cor}\label{thm3}
The lexicographic product $G\circ H$ has at least $|V(H)|$
arc-disjoint $S$-strong subgraphs for any $S\subseteq V(G\circ H)$
with $2\leq |S|\leq |V(G\circ H)|$, if one of the following
conditions holds:\\
$(i)$~$G$ is a strong symmetric digraph.\\
$(ii)$~$G$ is a strong semicomplete digraph and $H\not\cong
\overline{K_2}$.\\Moreover, these subgraphs can be found in
polynomial time.
\end{cor}


Recall that strong semicomplete compositions generalize strong quasi-transitive digraphs. Therefore, the following result holds by Theorem \ref{thmc}:
\begin{cor}\label{cor1}
Let $Q\not\in \mathcal{Q}_0$  
be a strong quasi-transitive digraph. We can in polynomial time find at least $n_0$ arc-disjoint $S$-strong subgraphs in $Q$ for
any $S\subseteq V(Q)$ with $2\leq |S|\leq |V(Q)|$.
\end{cor}

\section{Discussiones}

Let $G$ be a connected graph with $S\subseteq V(G)$. We say that a
set of edges $C$ of $G$ an {\em $S$-Steiner-cut} if there are at
least two components of $G\setminus C$ which contain vertices of
$S$. Similarly, let $D$ be a strong digraph and $S\subseteq V(D)$;
we say that a set of arcs $C$ of $D$ an {\em $S$-strong
subgraph-cut} if there are at least two strong components of
$D\setminus C$ which contain vertices of $S$.

Kriesell posed the following
well-known conjecture which concerns an approximate min-max relation
between the size of an $S$-Steiner-cut and the number of
edge-disjoint $S$-Steiner trees.

\begin{conj}\label{kriesell}\cite{Kriesell}
Let $G$ be a graph and $S\subseteq V(G)$ with $|S|\geq 2$. If every
$S$-Steiner-cut in $G$ has size at least $2{\ell}$, then $G$
contains $\ell$ pairwise edge-disjoint $S$-Steiner trees.
\end{conj}

Lau \cite{Lau} proved that the conjecture holds if every
$S$-Steiner-cut in $G$ has size at least $26{\ell}$. West and Wu \cite{West-Wu} improved the bound significantly by showing that the conjecture still holds if $26{\ell}$ is replaced by $6.5{\ell}$. So far the best bound $5{\ell}+4$ was obtained by DeVos, McDonald and Pivotto as follows.  

\begin{thm}\label{thm05}\cite{DeVos-McDonald-Pivotto}
Let $G$ be a graph and $S\subseteq V(G)$ with $|S|\geq 2$. If every
$S$-Steiner-cut in $G$ has size at least $5{\ell}+4$, then $G$
contains $\ell$ pairwise edge-disjoint $S$-Steiner trees.
\end{thm}







Similar to Theorem~\ref{thm05}, it is natural to study an approximate min-max relation
between the size of minimum $S$-strong subgraph-cut and the maximum
number of arc-disjoint $S$-strong subgraphs in a digraph $D$. Here is an interesting problem which is analogous to
Conjecture \ref{kriesell}.

\begin{op}\label{op1}
Let $D$ be a digraph and $S\subseteq V(D)$ with $|S|\geq 2$. Find a function
$f(\ell)$ such that the following holds: If every $S$-strong
subgraph-cut in $G$ has size at least $f(\ell)$, then $D$ contains
$\ell$ pairwise arc-disjoint $S$-strong subgraphs.
\end{op}

Note that there is a linear function $f(\ell)$ for a strong symmetric
digraph: Let $D=\overleftrightarrow{G}$ be a strong symmetric digraph and
$S\subseteq V(D)$. If every $S$-strong subgraph-cut in $D$ has size at least
$10{\ell}+8$, then $D$ contains $\ell$ pairwise arc-disjoint
$S$-strong subgraphs. The argument is as follows: Let $c_1$
and $c_2$ be the sizes of the minimum $S$-Steiner-cut in $G$ and the
minimum $S$-strong subgraph-cut in $D$, respectively. We deduce that
$c_1\geq \frac{c_2}{2}$. Indeed, let $C_1=\{e_i\mid 1\leq i\leq
c_1\}$ be the minimum $S$-Steiner-cut in $G$. Let $C'_1=\{a_i,
a'_i\mid 1\leq i\leq c_1\}$, where $a_i, a'_i$ be the two arcs in
$D$ corresponding to the edge $e_i$. It can be checked that $C'_1$
is an $S$-strong subgraph-cut of $D$. Hence, $c_2\leq |C'_1|= 2c_1$.
The assumption means that $c_2 \geq 10{\ell}+8$ and so $c_1\geq
\frac{c_2}{2} \geq 5{\ell}+4$. By Theorem \ref{thm05}, $G$ contains $\ell$ pairwise edge-disjoint $S$-Steiner trees. For
each $S$-Steiner tree, we can obtain an $S$-strong subgraph in $D$
by replacing each edge of this tree with the corresponding arcs of
both directions in $D$. Observe that we now obtain $\ell$
pairwise arc-disjoint $S$-strong subgraphs.

\vskip 1cm




\noindent{\large\bf Acknowledgements}

We thank Professor Gexin Yu and Professor Hehui Wu for their helpful suggestions and constructive comments. This work was supported by Zhejiang Provincial Natural Science Foundation of China (No. LY20A010013), National Natural Science Foundation of China (Nos.11871280 and 11971349) and Qing Lan Project.

\end{document}